\documentclass[11pt]{amsart}
\usepackage[T1]{fontenc}
\usepackage[margin=3cm]{geometry}
\usepackage{amsmath,amssymb,amsthm,mathtools}
\usepackage[hidelinks]{hyperref}
\usepackage{tikz}
\numberwithin{equation}{section}
\newtheorem{theorem}{Theorem}[section]
\newtheorem{proposition}[theorem]{Proposition}
\newtheorem{lemma}[theorem]{Lemma}
\newcommand{\R}{\mathbb R}
\title{Degenerate neckpinches at the first \texorpdfstring{$N$}{N} singular times}
\author{Jingwen Chen}
\address{Department of Mathematics, Faculty of Science,
University of Macau,\newline
E11, Avenida da Universidade, Taipa,
Macau}
\email{jingwenchen@um.edu.mo}
\date{}
\begin{document}
\begin{abstract}
For every integer $N\ge1$, we construct a smooth closed connected strictly mean-convex rotationally symmetric surface in $\R^3$ whose level set flow remains connected until extinction and has exactly one singular point at each of its first $N$ singular times. These singularities are degenerate neckpinches with the Angenent--Vel\'azquez asymptotics: the leading deviation from the cylinder under parabolic rescaling is a cubic Hermite term.
\end{abstract}
\subjclass[2020]{53E10}
\keywords{Mean curvature flow, degenerate neckpinches, rotationally symmetric surfaces}
\maketitle

\section{Introduction}\label{sec:introduction}

A smooth one-parameter family of surfaces $M(t)\subset\R^3$, $t\in(0,T)$, evolves by its mean curvature if
\[
 \partial_t x=\vec H(x),
\]
where $\vec H$ is the mean curvature vector. The monotonicity formula of Huisken \cite{Hui} relates singularities of this flow to self-similarly shrinking solutions. Colding--Minicozzi \cite{CM12} introduced entropy and showed that generalized cylinders are the only entropy stable singularity models. The work of Chodosh--Choi--Mantoulidis--Schulze \cite{CCMS24,CCMS26}, together with the multiplicity one theorem of Bamler--Kleiner \cite{BK}, shows that a mean curvature flow starting from a generic closed embedded surface in $\R^3$ develops only multiplicity-one spherical and cylindrical singularities. For such flows, Colding--Minicozzi \cite{CM16} proved smoothness at almost every time, and Choi--Haslhofer--Hershkovits \cite{CHH22} established mean convex neighborhoods and uniqueness of continuation through spherical and cylindrical singularities.

White \cite{Whi00,Whi03} proved that, for mean convex surfaces in $\R^3$, tangent flows at singular points are multiplicity-one round spheres and cylinders. A singularity with a round shrinking cylinder as a tangent flow is called a neckpinch. Colding--Minicozzi \cite{CM15} proved uniqueness of cylindrical tangent flows.

Altschuler--Angenent--Giga \cite{AAG} obtained a degenerate neckpinch using the interpolation method in their flow through singularities for surfaces of rotation. Angenent--Vel\'azquez \cite{AV} used the Wa\.{z}ewski box argument to construct closed examples of degenerate neckpinches with prescribed asymptotic expansions. In their examples where only one end disappears, parabolic rescalings at the singular point converge to a shrinking cylinder, whereas rescalings about the disappearing pole by its mean curvature converge to the bowl soliton. The remaining surface has a cusp at the singular time. Both works also discuss peanut solutions, in which the entire surface shrinks to a point without ever becoming convex. Angenent--Daskalopoulos--{\v S}e{\v s}um \cite{ADS} proved an instability result for a class of peanut solutions.

Ao Sun and Jinxin Xue \cite{SX} established normal-form asymptotics at cylindrical singularities and defined nondegeneracy through the leading asymptotic profile. They proved that nondegenerate singularities are isolated in space and have mean convex neighborhoods, and established their stability under small perturbations in \cite{SXII}. Under this nondegeneracy condition, Ao Sun, Zhihan Wang, and Jinxin Xue \cite{SWX} described the local geometry and topology of the flow through the singularity, showing that the local topological change agrees with surgery. In further work \cite{SWXreg}, they obtained detailed asymptotic profiles at cylindrical singularities and related them to the regularity of the cylindrical singular sets.

We call a singularity at $(p,T)$ locally Type~I if there are constants $r,C>0$ such that
\[
 |A(x,t)|^2\le\frac{C}{T-t},
 \qquad (x\in M(t)\cap B_r(p),\ T-r^2<t<T),
\]
and locally Type~II otherwise. Using a definition of nondegeneracy based on blow-up limits, Schulze--{\v S}e{\v s}um \cite{SS} proved that nondegenerate neckpinches are locally Type~I, whereas degenerate neckpinches are locally Type~II. At an isolated nondegenerate neckpinch, neither side of the collapsing neck disappears at the singular time. In contrast, at the degenerate neckpinches constructed here, one side of the neck collapses to the singular point.

Constructions prescribing the type or precise asymptotic behavior of singularities have largely focused on the first singular time. Tang-Kai Lee and Xinrui Zhao \cite{LZ} realized asymptotically conical shrinkers as singularity models of closed mean curvature flows. In subsequent joint work with them and Ao Sun \cite{CLSZ}, the author extended this result to a broader class of shrinkers and to flows in Riemannian manifolds. White \cite{Whi03,Whi15} established a general structure theory for subsequent singularities of mean convex flows. Much less is known about constructing flows with prescribed singularity type at a later singular time.

We study how to realize degenerate neckpinches with prescribed asymptotics at successive singular times of a mean curvature flow that remains connected through singularities. We use level set flow, a weak formulation of mean curvature flow that allows continuation through singularities.

\begin{theorem}\label{thm:main}
For every integer $N\ge1$, there is a smooth closed embedded connected strictly mean-convex rotationally symmetric surface $M(0)\subset\R^3$ whose level set flow has exactly one singular point at each of its first $N$ singular times and remains connected until extinction. These singularities are degenerate neckpinches with the $m=3$ asymptotics in \cite{AV}.
\end{theorem}

To the best of the author's knowledge, this is the first result prescribing a degenerate neckpinch at later singular times of a mean curvature flow from a smooth closed connected initial surface. Figure~\ref{fig:initial-surface} illustrates the mean curvature flow we construct for $N=2$, with first two singular times $T_1<T_2$. The drawings are schematic and not to scale, with the upper parts enlarged for clarity.
\begin{figure}[htbp]
\centering
\begin{tikzpicture}[yscale=1.08,line cap=round,line join=round,font=\small]
\path[fill=black!5,draw=black,line width=.75pt] (1.3,0.18) .. controls (1.89,0.18) and (2.34,0.52) .. (2.34,1.1) .. controls (2.34,1.78) and (1.54,1.83) .. (1.54,2.2) .. controls (1.54,2.5) and (1.99,2.63) .. (1.99,2.9) .. controls (1.99,3.23) and (1.42,3.31) .. (1.42,3.62) .. controls (1.42,3.85) and (1.68,3.85) .. (1.68,4.02) .. controls (1.68,4.26) and (1.5,4.43) .. (1.3,4.43) .. controls (1.1,4.43) and (0.92,4.26) .. (0.92,4.02) .. controls (0.92,3.85) and (1.18,3.85) .. (1.18,3.62) .. controls (1.18,3.31) and (0.61,3.23) .. (0.61,2.9) .. controls (0.61,2.63) and (1.06,2.5) .. (1.06,2.2) .. controls (1.06,1.83) and (0.26,1.78) .. (0.26,1.1) .. controls (0.26,0.52) and (0.71,0.18) .. (1.3,0.18) -- cycle;
\draw[black!40,densely dashed,line width=.3pt] (1.3,0.08) -- (1.3,4.67);
\draw[black!45,densely dashed,line width=.35pt] (2.34,1.1) arc[start angle=0,end angle=180,x radius=1.04cm,y radius=0.15cm];
\draw[black!45,line width=.35pt] (0.26,1.1) arc[start angle=180,end angle=360,x radius=1.04cm,y radius=0.15cm];
\draw[black!45,densely dashed,line width=.35pt] (1.99,2.9) arc[start angle=0,end angle=180,x radius=0.69cm,y radius=0.1cm];
\draw[black!45,line width=.35pt] (0.61,2.9) arc[start angle=180,end angle=360,x radius=0.69cm,y radius=0.1cm];
\draw[black!45,densely dashed,line width=.35pt] (1.68,4.02) arc[start angle=0,end angle=180,x radius=0.38cm,y radius=0.06cm];
\draw[black!45,line width=.35pt] (0.92,4.02) arc[start angle=180,end angle=360,x radius=0.38cm,y radius=0.06cm];
\node[anchor=north] at (1.3,-0.1) {$t=0$};
\path[fill=black!5,draw=black,line width=.75pt] (4.6,0.2) .. controls (5.175,0.2) and (5.615,0.535) .. (5.615,1.1) .. controls (5.615,1.76) and (4.83,1.86) .. (4.83,2.2) .. controls (4.83,2.5) and (5.26,2.61) .. (5.26,2.9) .. controls (5.26,3.21) and (4.655,3.34) .. (4.655,3.62) .. controls (4.655,3.77) and (4.82,3.77) .. (4.82,3.88) .. controls (4.82,4.02) and (4.715,4.12) .. (4.6,4.12) .. controls (4.485,4.12) and (4.38,4.02) .. (4.38,3.88) .. controls (4.38,3.77) and (4.545,3.77) .. (4.545,3.62) .. controls (4.545,3.34) and (3.94,3.21) .. (3.94,2.9) .. controls (3.94,2.61) and (4.37,2.5) .. (4.37,2.2) .. controls (4.37,1.86) and (3.585,1.76) .. (3.585,1.1) .. controls (3.585,0.535) and (4.025,0.2) .. (4.6,0.2) -- cycle;
\draw[black!40,densely dashed,line width=.3pt] (4.6,0.08) -- (4.6,4.67);
\draw[black!45,densely dashed,line width=.35pt] (5.615,1.1) arc[start angle=0,end angle=180,x radius=1.015cm,y radius=0.15cm];
\draw[black!45,line width=.35pt] (3.585,1.1) arc[start angle=180,end angle=360,x radius=1.015cm,y radius=0.15cm];
\draw[black!45,densely dashed,line width=.35pt] (5.26,2.9) arc[start angle=0,end angle=180,x radius=0.66cm,y radius=0.095cm];
\draw[black!45,line width=.35pt] (3.94,2.9) arc[start angle=180,end angle=360,x radius=0.66cm,y radius=0.095cm];
\draw[black!45,densely dashed,line width=.35pt] (4.82,3.88) arc[start angle=0,end angle=180,x radius=0.22cm,y radius=0.04cm];
\draw[black!45,line width=.35pt] (4.38,3.88) arc[start angle=180,end angle=360,x radius=0.22cm,y radius=0.04cm];
\node[anchor=north] at (4.6,-0.1) {$0<t<T_1$};
\path[fill=black!5,draw=black,line width=.75pt] (7.9,0.22) .. controls (8.46,0.22) and (8.89,0.55) .. (8.89,1.1) .. controls (8.89,1.74) and (8.12,1.88) .. (8.12,2.2) .. controls (8.12,2.5) and (8.52,2.6) .. (8.52,2.9) .. controls (8.52,3.05) and (8.22,3.16) .. (8.1,3.29) .. controls (7.9,3.51) and (7.9,3.62) .. (7.9,3.62) .. controls (7.9,3.62) and (7.9,3.51) .. (7.7,3.29) .. controls (7.58,3.16) and (7.28,3.05) .. (7.28,2.9) .. controls (7.28,2.6) and (7.68,2.5) .. (7.68,2.2) .. controls (7.68,1.88) and (6.91,1.74) .. (6.91,1.1) .. controls (6.91,0.55) and (7.34,0.22) .. (7.9,0.22) -- cycle;
\draw[black!40,densely dashed,line width=.3pt] (7.9,0.08) -- (7.9,4.67);
\draw[black!45,densely dashed,line width=.35pt] (8.89,1.1) arc[start angle=0,end angle=180,x radius=0.99cm,y radius=0.15cm];
\draw[black!45,line width=.35pt] (6.91,1.1) arc[start angle=180,end angle=360,x radius=0.99cm,y radius=0.15cm];
\draw[black!45,densely dashed,line width=.35pt] (8.52,2.9) arc[start angle=0,end angle=180,x radius=0.62cm,y radius=0.09cm];
\draw[black!45,line width=.35pt] (7.28,2.9) arc[start angle=180,end angle=360,x radius=0.62cm,y radius=0.09cm];
\node[anchor=north] at (7.9,-0.1) {$t=T_1$};
\path[fill=black!5,draw=black,line width=.75pt] (11.2,0.24) .. controls (11.745,0.24) and (12.165,0.565) .. (12.165,1.1) .. controls (12.165,1.72) and (11.35,1.93) .. (11.35,2.2) .. controls (11.35,2.46) and (11.68,2.53) .. (11.68,2.8) .. controls (11.68,3.05) and (11.42,3.31) .. (11.2,3.31) .. controls (10.98,3.31) and (10.72,3.05) .. (10.72,2.8) .. controls (10.72,2.53) and (11.05,2.46) .. (11.05,2.2) .. controls (11.05,1.93) and (10.235,1.72) .. (10.235,1.1) .. controls (10.235,0.565) and (10.655,0.24) .. (11.2,0.24) -- cycle;
\draw[black!40,densely dashed,line width=.3pt] (11.2,0.08) -- (11.2,4.67);
\draw[black!45,densely dashed,line width=.35pt] (12.165,1.1) arc[start angle=0,end angle=180,x radius=0.965cm,y radius=0.145cm];
\draw[black!45,line width=.35pt] (10.235,1.1) arc[start angle=180,end angle=360,x radius=0.965cm,y radius=0.145cm];
\draw[black!45,densely dashed,line width=.35pt] (11.68,2.8) arc[start angle=0,end angle=180,x radius=0.48cm,y radius=0.075cm];
\draw[black!45,line width=.35pt] (10.72,2.8) arc[start angle=180,end angle=360,x radius=0.48cm,y radius=0.075cm];
\node[anchor=north] at (11.2,-0.1) {$T_1<t<T_2$};
\end{tikzpicture}
\caption{}
\label{fig:initial-surface}
\end{figure}

\subsection{Proof ideas}

Starting with a sphere, we form successive families of initial surfaces by applying the local construction of Angenent--Vel\'azquez \cite{AV} near the upper pole, adding one neck and three parameters at each step. We choose the scales starting with the last singularity. Each modification also affects the later evolution. We therefore show that the modified flows are close in $C^\infty$ to the preceding flows at fixed regular times between the prescribed singular times, uniformly over the parameter families. For sufficiently small modifications, this preserves the strict Angenent--Vel\'azquez initial conditions for each later stage. We then use one Wa\.{z}ewski box argument, as in the author's joint work \cite{CLSZ}, to select all $3N$ parameters simultaneously and obtain the prescribed degenerate neckpinches.

The Sturmian theorem and barrier arguments exclude additional singular points or singular times before final extinction. Finally, the Angenent--Vel\'azquez estimates give the prescribed asymptotics at each neckpinch.

In Section~\ref{sec:preliminary}, we collect the results used below, including the Angenent--Vel\'azquez conditions and estimates. In Section~\ref{sec:nth}, we give the construction and prove Theorem~\ref{thm:main}.

\subsection*{Acknowledgment}
The author is supported by the University of Macau Research Grant (MYRG-\allowbreak SRG2026-\allowbreak 00039-\allowbreak FSC). The author thanks Ao Sun for his interest.

\subsection*{AI declaration}
The main methods in this paper are those used in the author's previous joint work \cite{CS23,CS24,CLSZ}: the use of the parabolic Sturmian theorem and comparison arguments for rotationally symmetric mean curvature flow, together with the Wa\.{z}ewski box argument. These methods were already used in \cite{AAG,AV}. Before using AI, the author had developed the main ideas and proof strategy and completed a draft. ChatGPT 6 Pro was used to assist with literature searches, streamlining the proof, and checking the arguments for possible gaps. The author checked and revised the arguments and takes full responsibility for the mathematical content and final exposition.
\par

\section{Preliminaries}\label{sec:preliminary}
\subsection{Rotationally symmetric mean curvature flow}\label{sec:rotational}

We choose the $y$-axis as the rotation axis, write a point of $\R^3$ as $(z,y)\in\R^2\times\R$, and set $x=|z|$. Following \cite{CS23,CS24}, we write the profiles as vertical graphs $x=u(y,t)$ and locally as horizontal graphs $y=v(x,t)$. When $M(t)$ is connected and smooth, we have
\[
 M(t)=\{(z,y): |z|=u(y,t),\ a(t)\le y\le b(t)\}.
\]
Here $a(t)$ and $b(t)$ are the axial coordinates of the lower and upper poles, and $u(a(t),t)=u(b(t),t)=0$. For surfaces, $C^k$ neighborhoods and $C^k$ convergence are understood in normal graph coordinates; for profiles, we use $u$ away from the poles and $v$ near smooth poles. We use the terminology concerning necks in \cite{AAG} and the definition of a degenerate neckpinch in \cite[Def.~1.1]{SS}.

By \cite{AAG,CS23}, the vertical graph function satisfies
\begin{equation}\label{eq:pde}
 u_t=\frac{u_{yy}}{1+u_y^2}-\frac1u.
\end{equation}
For the outward unit normal, the scalar mean curvature is
\begin{equation}\label{eq:H}
 H=\frac{1+u_y^2-u u_{yy}}{u(1+u_y^2)^{3/2}}.
\end{equation}
For the initial surfaces considered here, Altschuler--Angenent--Giga \cite{AAG} proved that the level set flow is nonfattening and its continuation through singularities is unique in the level set sense. Their analysis also shows that the flow remains smooth away from the rotation axis, has finitely many singular times, and becomes singular by collapsed necks or extinction of a whole component. A component with no neck can become singular only by extinction at a single point. They also proved that the number of necks does not increase. If the region enclosed by a component has nonempty interior at a singular time, its surviving components are smooth immediately afterwards, and the total number of necks decreases by at least the number of distinct singular points on that component. On fixed axial intervals, we apply the parabolic Sturmian theorem as in \cite{Ang88,CS23}.

\subsection{Rescaling and the Angenent--Vel\'azquez estimates}\label{sec:AV-estimates}
Fix a prescribed singular time $T>0$ and an axial coordinate $Y$ centered at the prescribed singular point, with the upper pole in $Y<0$. Following Angenent--Vel\'azquez \cite{AV}, set
\[
 \tau=-\log(T-t),\qquad
 \rho=\frac{Y}{\sqrt{T-t}},\qquad
 \Phi(\rho,\tau)=\frac{u(Y,t)}{\sqrt{T-t}}.
\]
Our coordinates $Y,\rho$ correspond to their $x,y$, respectively. Set $\psi=\Phi-\sqrt2$. The linearized operator at the cylinder $\Phi\equiv\sqrt2$ has Hermite eigenfunctions:
\[
 \mathcal A=\partial_{\rho\rho}-\frac\rho2\partial_\rho+1,
 \qquad
 \mathcal A\operatorname{Hm}_j=\left(1-\frac j2\right)\operatorname{Hm}_j.
\]
The polynomials $\operatorname{Hm}_j$ have leading coefficient one and are orthogonal in $L^2(\R,e^{-\rho^2/4}\,d\rho)$, whose inner product and norm we use below. The first four are
\[
 \operatorname{Hm}_0=1,\quad
 \operatorname{Hm}_1=\rho,\quad
 \operatorname{Hm}_2=\rho^2-2,\quad
 \operatorname{Hm}_3=\rho^3-6\rho.
\]
We use the case $m=3$ of \cite{AV}: the cubic coefficient in the initial polynomial is prescribed, while three parameters vary the coefficients of $\operatorname{Hm}_0,\operatorname{Hm}_1,\operatorname{Hm}_2$.

Fix a closed rotationally symmetric surface $\Gamma$, smooth except at its pole $Y=0$, with local profile $u=\sqrt{K_0}Y^{3/2}$ for $Y\ge0$, where $K_0>0$. Choose $\delta,R>0$ sufficiently small and an even cutoff function $\eta_0\in C_c^\infty(\R)$ with $0\le\eta_0\le1$, equal to one on $[-1,1]$ and zero outside $[-2,2]$. Set
\[
 \xi=\frac{Y}{(T-t)^{1/3}}=\rho e^{-\tau/6},\qquad
 \eta(\rho,\tau)=\eta_0\!\left(\frac\xi R\right).
\]
We extend $\eta\psi$ by zero outside the profile interval. For $M_0>0$, define the normalized coefficient vector
\begin{equation}\label{eq:AV-low}
 V(t)=\frac1{M_0(T-t)}
 \left(
 \frac{\langle\eta\psi,\operatorname{Hm}_0\rangle}{\|\operatorname{Hm}_0\|},
 \frac{\langle\eta\psi,\operatorname{Hm}_1\rangle}{\|\operatorname{Hm}_1\|},
 \frac{\langle\eta\psi,\operatorname{Hm}_2\rangle}{\|\operatorname{Hm}_2\|}
 \right)\in\R^3.
\end{equation}
The additional bound $|V(t)|\le1$ is the coefficient condition in \cite[(3.1)]{AV}, with our cutoff.

We call the following requirements \textbf{Condition~(A)}: at a starting time $t_0<T$, the surface satisfies the initial conditions in \cite[(3.2)--(3.9)]{AV} for $T,\Gamma$ and fixed constants, using $\eta$ in all Hermite projections, and $u(Y,t_0)$ is defined for $|Y|\le\tfrac12(2(T-t_0)/K_0)^{1/3}$.

For the choices of $R$ below, $\eta$ vanishes near the profile endpoints, so integration by parts produces no boundary terms. Its derivatives are supported where $Re^{\tau/6}\le|\rho|\le2Re^{\tau/6}$. The pointwise bounds in \cite[Lem.~4.1]{AV} give a Gaussian $L^2$ cutoff error $O(\exp(-R^2e^{\tau/3}/16))$. With $\Phi\ge1$ on $|\xi|\le3R$, the projection calculation in that lemma applies with our cutoff. Taking $M_0$ and the initial rescaled time sufficiently large gives, at $|V|=1$,
\begin{equation}\label{eq:AV-outward}
 \frac{d}{d\tau}|V|\ge\frac34,
\end{equation}
including at the initial time. Projection onto the cubic mode and the modes $j\ge4$ gives the corresponding estimates from the same lemma with our cutoff.

The inner estimates and barriers in \cite[Lem.~4.2, 6.1--6.2]{AV} preserve $\Phi\ge1$ on $|\xi|\le3R$ when $2-(K_0+3\delta)(3R)^3>1$. While $|V|\le1$, the argument in \cite[p.~47]{AV} gives a curvature bound outside a fixed neighborhood of the cusp for sufficiently small $T-t_0$. This bound and the required smallness are uniform for fixed barrier constants and common initial $C^3$ and curvature bounds there. Thus the continuation and exit conclusions in \cite[Lem.~3.1]{AV} hold with our cutoff: for initial data satisfying Condition~(A) strictly and $|V(t_0)|\le1$, with $T-t_0$ sufficiently small, the flow either is classical on $[t_0,T)$ with $|V|\le1$ throughout, or first reaches $|V|=1$ at $t\in[t_0,T)$, where it extends smoothly and crosses strictly outward by \eqref{eq:AV-outward}.

\subsection{Smooth convergence at positive times}
Let $M_{\mathbf p}(t)$, with $\mathbf p$ in a compact parameter set, be a smooth family of classical flows on $[0,T]$ of closed connected strictly mean-convex rotationally symmetric surfaces with vertical graph profiles, and let $M^\sigma_{\mathbf p}(t)$ be the level set flows from smooth families of such initial surfaces.

\begin{lemma}\label{lem:regularity}
Suppose $M^\sigma_{\mathbf p}(0)$ agrees with $M_{\mathbf p}(0)$ outside a ball of radius $\sigma$ centered on $M_{\mathbf p}(0)$, and the total added and removed area tends uniformly to zero as $\sigma\downarrow0$. For every sufficiently small fixed $s>0$, $M^\sigma_{\mathbf p}(t)\to M_{\mathbf p}(t)$ in $C^\infty$ as single normal graphs, uniformly in $\mathbf p$ and $t\in[s,2s]$, as $\sigma\downarrow0$. For fixed sufficiently small $\sigma$, the slices $M^\sigma_{\mathbf p}(t)$ depend continuously on $\mathbf p$ in $C^\infty$.
\end{lemma}

\begin{proof}
The enclosed regions differ only inside the modification ball. Comparison with the flows starting from the inward and outward parallel surfaces of $M_{\mathbf p}(0)$ at distance $2\sigma$, together with classical continuous dependence, gives Hausdorff convergence on compact time intervals in $[0,T)$. Choose $\gamma>0$ below the threshold in White's local regularity theorem \cite{Whi05} and below $4/e-1$. Smoothness and compactness bound the initial Gaussian integrals of $M_{\mathbf p}(0)$ by $1+\gamma/4$, uniformly over all centers and parameters at sufficiently small scales. Fix a sufficiently small $s>0$. For backward time parameters in $[s/2,3s]$, the Gaussian kernel is bounded by $(2\pi s)^{-1}$, so the area hypothesis bounds the corresponding integrals of $M^\sigma_{\mathbf p}(0)$ by $1+\gamma/2$. By \cite{Hui,Whi00}, monotonicity remains valid through singularities; together with the singular density lower bound $4/e$ in \cite{Whi03}, it excludes singularities on $[s/2,3s]$. White's curvature estimate and the interior derivative estimates in \cite{EH91} then give uniform bounds for all derivatives of curvature on $[s,2s]$.

The vertical graph property excludes multiple sheets away from the poles, while the curvature bounds exclude additional components near the poles. Thus the convergence is as a single normal graph. For fixed $\sigma$ and a convergent parameter sequence, comparison with the inward and outward parallel surfaces of the limiting initial surface $M^\sigma_{\mathbf p}(0)$, together with nonfattening and the monotone convergence in \cite{AAG}, gives Hausdorff convergence on $[s,2s]$. The uniform derivative bounds then give $C^\infty$ continuity of the slices.
\end{proof}

\section{The construction}\label{sec:nth}

We first prove a local modification lemma that ensures Condition~(A) and the support and area bounds required by Lemma~\ref{lem:regularity}. We then construct the initial families and select the parameters by a Wa\.{z}ewski box argument. Finally, we exclude additional singularities and derive the asymptotics to complete the proof of Theorem~\ref{thm:main}.

\subsection{The local modification}\label{sec:asymptotics}

See Figure~\ref{fig:local-modification} for an illustration of the local modification.

Let $N_\theta$, with $\theta$ in a compact parameter set, be a smooth family of closed strictly mean-convex rotationally symmetric surfaces with vertical graph profiles. Suppose all surfaces $N_\theta$ coincide near the upper pole and are strictly convex there. Let $B^3\subset\R^3$ denote the closed unit ball.
We construct $\Gamma_\theta$ from $N_\theta$ by replacing the upper cap with one that is smooth except at its pole and has the common cusp profile described in Subsection~\ref{sec:AV-estimates}.

\begin{lemma}\label{lem:AV-local}
For every sufficiently small $\sigma>0$ and $T>0$ sufficiently small depending on $\sigma$, we can replace an upper cap of each $N_\theta$ to obtain a smooth family of closed rotationally symmetric surfaces $N_{\mathbf p,\theta}$, $\mathbf p\in B^3$, with $H>0$ and vertical graph profiles. They agree with $N_\theta$ outside a ball of radius $C\sigma$ centered at the upper pole of $N_\theta$, with total added and removed area at most $C\sigma^2$. Here $C$ is independent of $\sigma,T,\mathbf p,\theta$. The profile of each new cap has exactly one positive local minimum and one local maximum, and all surfaces $N_{\mathbf p,\theta}$ share a strictly convex cap near the new upper pole.

These surfaces satisfy Condition~(A) strictly and uniformly at time $0$ for $T$ and $\Gamma_\theta$. The axial coordinate, cusp profile and constants are independent of $\mathbf p,\theta$, and
\begin{equation}\label{eq:AV-first-degree}
 \sup_{\mathbf p,\theta}|V(\mathbf p,\theta,0)-2\mathbf p|<\tfrac14.
\end{equation}
\end{lemma}
\begin{proof}
Choose a sufficiently small circle of radius less than $\sigma/4$ in the common upper cap of $N_\theta$, where $v_x(x,0)<0$ and $v_{xx}(x,0)<0$. We first construct $\Gamma_\theta$ by solving $u_y=-q(u)$, matching the original profile at this circle. Here $q(x)>0$ is the slope magnitude at radius $x$, chosen so that
\[
 q(x)=
 \begin{cases}
 \tfrac32K_0^{1/3}x^{1/3},&\text{near the axis},\\
 -1/v_x(x,0),&\text{near the chosen circle}.
 \end{cases}
\]
For large $K_0$, choose $\log(1+q^2)$ as a fixed smooth minimum of the logarithms of one plus the squares of the two displayed slopes. Both logarithms have derivative at most $2/3$ with respect to $\log x$. The derivative of their smooth minimum is a convex combination of these derivatives, giving $xqq'\le\tfrac13(1+q^2)$. Since $-1/v_x(x,0)\asymp x^{-1}$, the transition occurs at $x\asymp K_0^{-1/4}$; choose $K_0$ large enough that it lies inside half the chosen circle's radius. Hence
\[
 1+u_y^2-u u_{yy}=1+q^2-xqq'
 \ge\tfrac23(1+q^2)>0.
\]
Set $Y=0$ at the new endpoint and orient $Y$ toward the unmodified surface. Rotating this profile gives $\Gamma_\theta$, which agrees with $N_\theta$ outside the modified cap and has profile $u=\sqrt{K_0}Y^{3/2}$ near $Y=0$. The transition and cusp have axial length $O(K_0^{-1/2})$, whereas the part of the original cap between the transition and the chosen circle has axial length bounded by a constant times the square of the chosen circle's radius. Thus the modified part of $\Gamma_\theta$ has radial size at most $\sigma/4$ and axial extent $C\sigma^2$, with $C$ independent of $\sigma,K_0,\theta$.

Choose the constants in the order used in the proof of \cite[Lem.~3.1]{AV}, with the cutoff in Subsection~\ref{sec:AV-estimates}. The radius $\varepsilon>0$ of the neighborhood of $\Gamma_\theta$ required in Condition~(A) will be chosen below. In choosing the constants, require $\delta<K_0/10$, $R$ sufficiently small that the sections $Y=4R$ and $Y=2(R^2/(K_0-3\delta))^{1/3}$ lie in the exact cusp region, and $10R<\tfrac12(2/K_0)^{1/3}$.

First smooth the cusp of $\Gamma_\theta$ using a profile $f$ equal to $\sqrt{2T+K_0Y^3}$ for $-(2T/K_0)^{1/3}\le Y\le2R$ and to the profile of $\Gamma_\theta$ for $Y\ge3R$, joining them with a fixed smooth cutoff on $2R<Y<3R$. On this annulus, the error relative to $\Gamma_\theta$ is $O(T)$ in every $C^k$ norm. Near the new upper pole,
\[
 Y(x)=\left(\frac{x^2-2T}{K_0}\right)^{1/3}
\]
is smooth, even, and strictly convex. This cap is independent of $\mathbf p$ and $\theta$.

To construct the profile $u$ of $N_{\mathbf p,\theta}$, modify $f$ near $Y=0$ as in \cite[Sec.~3]{AV}. For $\mathbf p=(p_1,p_2,p_3)\in B^3$, set $\rho=Y/\sqrt T$, $\xi=Y/T^{1/3}$ and prescribe $u/\sqrt T$ on $|\xi|\le5R$ by
\[
 \sqrt2+\frac{K_0}{2\sqrt2}\sqrt T\operatorname{Hm}_3(\rho)
       +2M_0T\sum_{j=0}^2p_{j+1}\frac{\operatorname{Hm}_j(\rho)}{\|\operatorname{Hm}_j\|}.
\]
Extend $u/\sqrt T$ as the weighted sum of the displayed polynomial and $f/\sqrt T$, with respective weights $\eta_0(\xi/(5R))$ and $1-\eta_0(\xi/(5R))$. For small $T$, the profile domain strictly contains the interval in Condition~(A).

\begin{figure}[htbp]
\centering
\begin{tikzpicture}[line cap=round,line join=round,font=\small]
\node at (.95,3.8) {$N_\theta$};
\node at (4.15,3.8) {$\Gamma_\theta$};
\node at (7.35,3.8) {profile $f$};
\node at (10.55,3.8) {$N_{\mathbf p,\theta}$};
\draw[->,line width=.45pt] (1.6,3.8) -- (3.5,3.8)
 node[midway,above=3pt,font=\scriptsize] {form cusp};
\draw[->,line width=.45pt] (4.8,3.8) -- (6.6,3.8)
 node[midway,above=3pt,font=\scriptsize] {round pole};
\draw[->,line width=.45pt] (8.1,3.8) -- (9.8,3.8)
 node[midway,above=3pt,font=\scriptsize,align=center] {Hermite\\perturbation};
\begin{scope}[shift={(.95,0)}]
 \draw[black!40,densely dashed,line width=.3pt] (0,0) -- (0,3.4);
 \foreach \s in {-1,1} {
  \begin{scope}[xscale=\s]
   \draw[line width=.7pt] (0,1.52) .. controls (.2,1.52) and (.4,1.38) .. (.6,1.1)
    .. controls (.75,.89) and (.86,.51) .. (.93,.12);
  \end{scope}
 }
 \draw[black!50,densely dotted,line width=.4pt] (-.6,1.1) -- (.6,1.1);
 \draw[black!50,line width=.3pt] (.6,1.1) -- (.6,.83);
 \draw[<->,black!60,line width=.35pt] (0,.9) -- (.6,.9)
  node[midway,below=2pt,fill=white,inner sep=1pt,font=\scriptsize] {$<\sigma/4$};
 \node at (0,-.25) {original cap};
\end{scope}
\begin{scope}[shift={(4.15,0)}]
 \draw[black!40,densely dashed,line width=.3pt,->] (0,3.4) -- (0,0)
  node[pos=.85,right,font=\scriptsize] {$Y$};
 \draw[black!35,densely dotted,line width=.35pt] (-.6,2.45) -- (.6,2.45);
 \node[anchor=east,font=\scriptsize] at (-.65,2.45) {$Y=0$};
 \foreach \s in {-1,1} {
  \begin{scope}[xscale=\s]
   \draw[line width=.7pt] (0,2.45) .. controls (0,1.95) and (.45,1.31) .. (.6,1.1)
    .. controls (.75,.89) and (.86,.51) .. (.93,.12);
  \end{scope}
 }
 \fill (0,2.45) circle[radius=1.1pt];
 \node at (0,-.25) {cusp};
\end{scope}
\begin{scope}[shift={(7.35,0)}]
 \draw[black!40,densely dashed,line width=.3pt,->] (0,3.4) -- (0,0)
  node[pos=.85,right,font=\scriptsize] {$Y$};
 \draw[black!35,densely dotted,line width=.35pt] (-.6,2.45) -- (.6,2.45);
 \node[anchor=east,font=\scriptsize] at (-.65,2.45) {$Y=0$};
 \foreach \s in {-1,1} {
  \begin{scope}[xscale=\s]
   \draw[line width=.7pt] (0,3.25) .. controls (.15,3.25) and (.24,3.17) .. (.28,3.08)
    .. controls (.32,2.99) and (.36,2.8) .. (.36,2.45)
    .. controls (.36,1.75) and (.45,1.31) .. (.6,1.1)
    .. controls (.75,.89) and (.86,.51) .. (.93,.12);
  \end{scope}
 }
 \node at (0,-.25) {$x=f(Y)$};
\end{scope}
\begin{scope}[shift={(10.55,0)}]
 \draw[black!40,densely dashed,line width=.3pt,->] (0,3.4) -- (0,0)
  node[pos=.85,right,font=\scriptsize] {$Y$};
 \draw[black!35,densely dotted,line width=.35pt] (-.6,2.45) -- (.6,2.45);
 \node[anchor=east,font=\scriptsize] at (-.65,2.45) {$Y=0$};
 \foreach \s in {-1,1} {
  \begin{scope}[xscale=\s]
   \draw[black!45,densely dashed,line width=.6pt] (0,3.25)
    .. controls (.15,3.25) and (.24,3.17) .. (.28,3.08)
    .. controls (.32,2.99) and (.36,2.8) .. (.36,2.45)
    .. controls (.36,1.75) and (.45,1.31) .. (.6,1.1);
   \draw[line width=.7pt] (0,3.25) .. controls (.15,3.25) and (.24,3.17) .. (.28,3.08)
    .. controls (.32,2.99) and (.43,2.93) .. (.43,2.8)
    .. controls (.43,2.6) and (.23,2.25) .. (.23,2)
    .. controls (.23,1.65) and (.45,1.31) .. (.6,1.1)
    .. controls (.75,.89) and (.86,.51) .. (.93,.12);
  \end{scope}
 }
 \node[anchor=east,font=\scriptsize] at (.3,2.8) {$f$};
 \draw[black!65,line width=.35pt] (.29,2) -- (.83,2.12)
  node[anchor=west,font=\scriptsize] {neck};
 \node at (0,-.25) {$x=u(Y)$};
\end{scope}
\end{tikzpicture}
\caption{Upper parts of the profiles in the local modification. The drawings are schematic and not to scale.}
\label{fig:local-modification}
\end{figure}
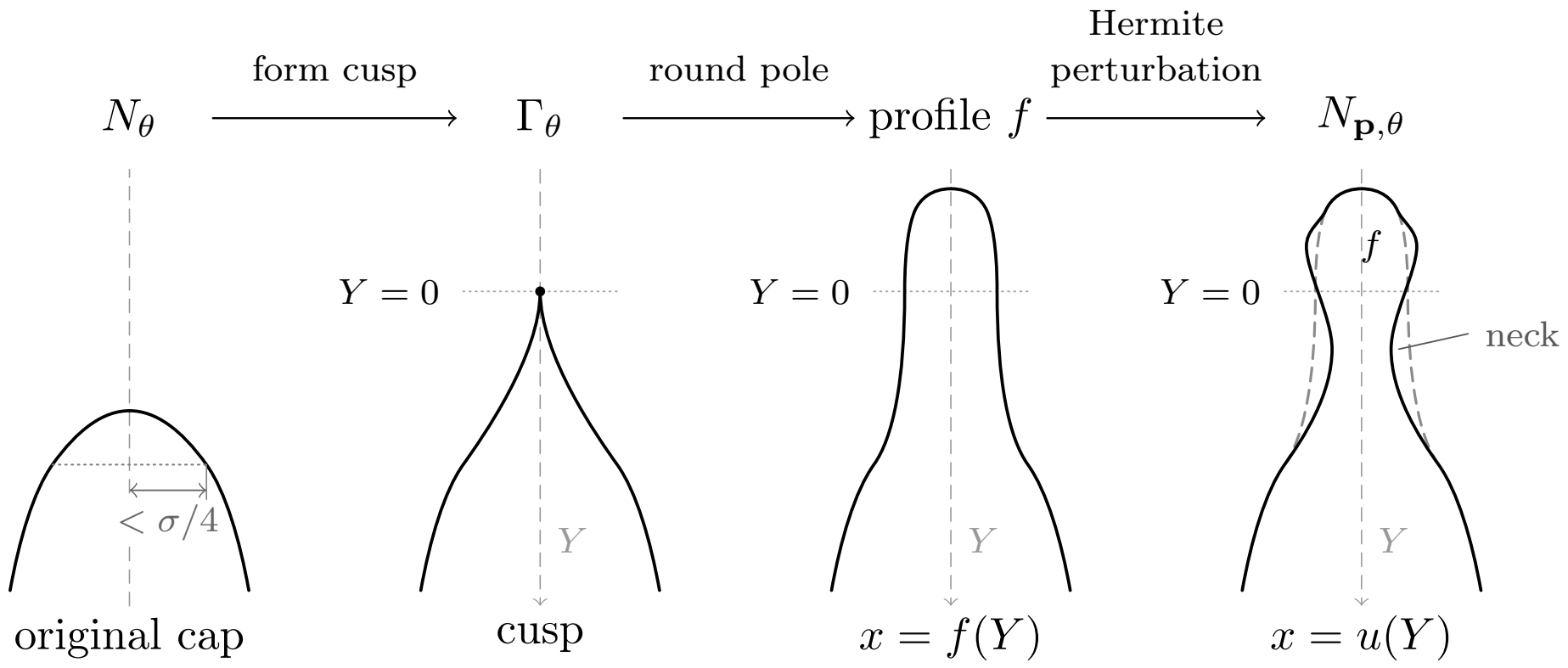

The explicit profiles and the outer interpolation give the neighborhood condition in (A):
\[
 \sup_{X\in N_{\mathbf p,\theta}}
       \operatorname{dist}(X,\Gamma_\theta)\le CT^{1/3}<\varepsilon
\]
for sufficiently small $T$, uniformly in $\mathbf p,\theta$.

The inner bounds and square-root comparisons away from the pole follow as in \cite[pp.~31--32]{AV}. The Hermite projection calculation in \cite[Prop.~3.4]{AV}, with the Gaussian cutoff error in Subsection~\ref{sec:AV-estimates}, gives $V=2\mathbf p+o(1)$, the $\operatorname{Hm}_3$ coefficient $\frac{\sqrt T}{2\sqrt2}(K_0+o(1))$ in $\eta\psi$, and the strict bounds on the higher Hermite modes, uniformly in $\mathbf p,\theta$. In particular, \eqref{eq:AV-first-degree} holds.

On $\xi\le-10R$, $u=f$, and the pole comparisons are those for the model inverse profile in \cite[Sec.~5]{AV}. Replacing $K_0$ by $K_0-\delta$ and $K_0+\delta$ gives strict lower and upper graphs, respectively, including at the axis, since $u/\sqrt T$ stays bounded away from $\sqrt2$ on this cap. The overlap comparisons follow as in \cite[pp.~31--32]{AV}.

Direct differentiation gives the required derivative signs on $|\xi|\le5R$ and $5R\le\xi\le10R$. On $-10R\le\xi\le-5R$, Taylor expansion shows that the first two derivatives of $u/\sqrt T$ with respect to $\xi$ differ from those of $\sqrt{2+K_0\xi^3}$ by $O_{K_0}(R^5)+o(1)$ and $O_{K_0}(R^4)+o(1)$, respectively, as $T\downarrow0$, uniformly in the parameters. The model's first derivative is at least $K_0R^2$ and its second derivative is at most $-K_0R$, so both signs persist for $R$ small and then $T$ small. For $\xi\le-10R$, we have $u_Y>0$ and $u_{YY}<0$; on $2R<Y<3R$, monotonicity persists by $C^\infty$ convergence to the cusp. Thus the monotonicity and concavity requirements in Condition~(A) hold. The two inner critical points converge to $\rho=\pm\sqrt2$, giving exactly one minimum and one maximum in the new cap.

To verify $H>0$ by \eqref{eq:H}, direct calculation gives $1+f'^2-ff''\ge1-C(K_0^2T)^{1/3}$ wherever $f^2=2T+K_0Y^3$, and $|u u_{YY}|=O(T^{1/3})$ on $|\xi|\le10R$. On $2R<Y<3R$, positivity follows from $C^2$ convergence to $\Gamma_\theta$.

For small $T$ depending on $\sigma$, the neighborhood estimate above and the size bounds for $\Gamma_\theta$ place the modification in a ball of radius $C\sigma$ centered at the original upper pole, with $C$ independent of $\sigma,T,\mathbf p,\theta$. The profile of the new cap has at most three monotone pieces. Integrating over this cap gives
\[
 2\pi\int u\sqrt{1+u_y^2}\,dy
 \le2\pi\int u\,dy+2\pi\int u\,|du|\le C\sigma^2.
\]
The area of the removed cap is also at most $C\sigma^2$.

The boundary circles at $x=R$ of the barriers in \cite[Sec.~5]{AV} have limits in the common cusp region, disjoint from every $\Gamma_\theta$; see \cite[proof of Lem.~6.1]{AV}. Choose $\varepsilon$ smaller than one quarter of their common positive distance from these surfaces. For small $T$, the circles remain outside the $2\varepsilon$-neighborhood of every $\Gamma_\theta$ throughout $0\le t<T$. The $O(T)$ error on the remaining fixed comparison interval preserves the strict initial barrier inequalities.

The common cusp and cutoff make the constants in \cite[Lem.~6.1--6.2]{AV} independent of $\theta$; impose also $T\le\varepsilon^2/4$. Outside a fixed neighborhood of the cusp, fix a common $C^3$ neighborhood of the portions of $\Gamma_\theta$ with a common curvature bound. The uniform $C^\infty$ convergence of $N_{\mathbf p,\theta}$ to $\Gamma_\theta$ there as $T\downarrow0$ makes the exterior curvature estimate in Subsection~\ref{sec:AV-estimates} applicable for small $T$. All requirements hold for $T$ sufficiently small depending on $\sigma$, uniformly in $\mathbf p,\theta$. The constants in Condition~(A) may depend on $\sigma$, whereas $C$ in the support and area bounds does not.
\end{proof}

\subsection{The initial surfaces}\label{sec:initial-surfaces}

Fix $N\ge1$, write $\mathbf p=(\mathbf p_1,\ldots,\mathbf p_N)\in(B^3)^N$, and set $s_0=0$. Starting from a round sphere $F_{N+1}$ centered at the origin, we apply Lemma~\ref{lem:AV-local} near the upper pole for $i=N,\ldots,1$ to construct the initial families $F_i(\mathbf p_i,\ldots,\mathbf p_N)$. Write $F_i(t)$ for their level set flows. Each application leaves the lower spherical cap unchanged and produces a compact smooth family with at most $N-i+1$ necks and a common strictly convex cap near the upper pole.

For $i<N$, suppose we have completed step $i+1$. Then $F_{i+1}$, $T_{i+1},s_{i+1}$ and the constants used at that step in Condition~(A) and the estimates of \cite{AV} are fixed. At each step, we choose $s_i$ first and then $T_i$, independently of $\mathbf p$, with
\[
 0=s_0<T_1<s_1<\cdots<T_N<s_N.
\]
Condition~(A) at this step uses $T_i$ and the axial coordinate, surfaces $\Gamma_\theta$ and constants chosen in that application of Lemma~\ref{lem:AV-local}; write $V_i$ for the corresponding coefficient vector.

Fix once and for all a section $y=a_-$ near the lower pole in the unchanged lower spherical cap. At step $i$, choose $y=a_i$ in the common upper cap of $F_{i+1}$, above all critical points of its profiles. Choose $a_i$ sufficiently close to the upper pole that the profile is strictly decreasing between the pole and its reflection across $y=a_i$. For $i<N$, require also $a_i>a_{i+1}$.

To keep the profile positive on $[a_-,a_i]$, choose finitely many closed balls centered on the axis, strictly inside every region enclosed by $F_{i+1}$, whose concentric open balls of half the radii cover this interval. Such balls exist by compactness and the uniform positive lower bound for the profiles. Include a ball centered at the origin to prevent extinction before $2s_i$. All modifications from step $i$ onwards avoid these balls.

Choose $s_i>0$ sufficiently small so that the flows $F_{i+1}(t)$ are smooth for $0\le t\le3s_i$, with profiles strictly decreasing above $a_i$ and strictly convex caps near both poles. Require also that Lemma~\ref{lem:regularity} applies at time $s_i$ and that each chosen ball, evolved by mean curvature flow, retains more than half its initial radius for $0\le t\le2s_i$. When $i<N$, also impose $s_i<\min\{T_{i+1}/3,s_{i+1}/2\}$ and use classical continuous dependence and compactness to ensure that $F_{i+1}(s_i)$ satisfies Condition~(A) strictly for step $i+1$ and
\[
 \sup_{\mathbf p_{i+1},\ldots,\mathbf p_N}
 |V_{i+1}(s_i)-2\mathbf p_{i+1}|<\tfrac13.
\]
Here we use the constants and rescaling fixed at step $i+1$. Condition~(A) allows the coefficient of $\operatorname{Hm}_3$ in $\eta\psi$, multiplied by $2\sqrt2e^{\tau/2}$, to vary in $[K_0-\delta/2,K_0+\delta/2]$. The pole comparisons persist in inverse graph coordinates, as in the proof of Lemma~\ref{lem:AV-local}. Fix a common $C^3$ neighborhood of these slices in which Condition~(A) remains strict, $|V_{i+1}-2\mathbf p_{i+1}|<\tfrac12$, and the portions outside the fixed cusp neighborhood stay in the $C^3$ neighborhood chosen at step $i+1$.

For every $i$, choose common $C^3$ neighborhoods of $F_{i+1}(t)$, $t\in[s_i,2s_i]$, preserving strict decrease above $a_i$ and strict convexity near the upper pole. By Lemma~\ref{lem:regularity}, fix upper bounds on the modification radius and the total added and removed area that preserve these profile properties on $[s_i,2s_i]$ and, for $i<N$, place the modified slice at $s_i$ in the chosen neighborhood of $F_{i+1}(s_i)$, uniformly in the parameters.

Apply Lemma~\ref{lem:AV-local} to $F_{i+1}$ with new parameter $\mathbf p_i$ and $\theta=(\mathbf p_{i+1},\ldots,\mathbf p_N)$. Choose the modification radius sufficiently small so that the modified region lies above $a_i$, is disjoint from all balls chosen so far, and preserves the comparisons in \eqref{eq:reflection}. Require its diameter and the total added and removed area to be less than half the respective bounds just fixed, and less than $1/(2N)$ of each corresponding bound fixed at a larger index. For this radius, fix the surfaces $\Gamma_\theta$, axial coordinate and Angenent--Vel\'azquez constants as in the proof of the lemma, together with a common $C^3$ neighborhood of their portions outside a fixed cusp neighborhood, with a common curvature bound. Finally choose $T_i<s_i/10$ sufficiently small and let $F_i$ be the resulting family. Let $q_i>a_i$ be the axial coordinate of the prescribed singular point, so $Y=q_i-y$. The surfaces coincide near the new upper pole, where they are strictly convex.

Set $M_{\mathbf p}(0)=F_1(\mathbf p)$ and denote its level set flow by $M_{\mathbf p}(t)$. For each $i$, regard $F_{i+1}$ as a family on $(B^3)^N$ independent of $\mathbf p_1,\ldots,\mathbf p_i$. The changes from $F_{i+1}$ to $F_1$ occur at steps $i,i-1,\ldots,1$. Since each modification displaces the upper pole by at most its support diameter, their union lies in a ball centered at the upper pole of $F_{i+1}$ with radius at most the sum of these diameters. This sum and the total added and removed area remain below the bounds fixed at step $i$. By Lemma~\ref{lem:regularity}, $M_{\mathbf p}(t)$ is smooth and connected on $[s_i,2s_i]$, with profiles strictly decreasing above $a_i$ and caps strictly convex near the upper pole. These slices depend continuously on $\mathbf p$ in $C^\infty$. For $i\ge2$, the neighborhood chosen at step $i-1$ ensures that $M_{\mathbf p}(s_{i-1})$ satisfies Condition~(A) strictly for step $i$ and
\begin{equation}\label{eq:AV-second-degree}
 \sup_{\mathbf p}|V_i(s_{i-1})-2\mathbf p_i|<\tfrac12.
\end{equation}
For $i=1$, the corresponding statements at $s_0=0$ follow from Lemma~\ref{lem:AV-local}.

Write $u(y,t)$ for the profile on each connected component of $M_{\mathbf p}(t)$. Initially,
\begin{equation}\label{eq:reflection}
 \begin{aligned}
 u(y,0)&<u(2a_i-y,0)&& (a_i<y<b(0)),\\
 u(y,0)&<u(2a_--y,0)&& (a(0)<y<a_-).
 \end{aligned}
\end{equation}
For the unmodified caps these inequalities follow from strict monotonicity. Each sufficiently small upper modification preserves the comparisons fixed at earlier steps: its profile radii tend uniformly to zero with the modification radius, while the reflected profile has a uniform positive lower bound at the corresponding points. The reflected images of the two cap intervals lie in $[a_-,a_i]$, where the balls chosen at step $i$ keep the profile positive for $0\le t\le2s_i$. Comparison on the positive profile intervals, including after singularities above $a_i$, and the Hopf lemma give
\begin{equation}\label{eq:section-sign}
 u_y(a_i,t)<0,\qquad u_y(a_-,t)>0
 \qquad(0\le t\le2s_i).
\end{equation}

Below $a_-$, the Sturmian theorem and \eqref{eq:section-sign} keep the profile increasing. Since $u(a_-,t)$ stays uniformly positive, the lower cap is a horizontal graph over a fixed disk. Its height is bounded by the initial enclosed regions, so the graphical interior estimates in \cite{AAG,EH91} bound curvature near the lower pole. On $[a_-,a_i]$, the shrinking balls keep $u$ uniformly positive; this part of the surface stays away from the rotation axis and is smooth by Subsection~\ref{sec:rotational}. Thus no singularity occurs at or below $a_i$ for $0\le t\le2s_i$.

At time $s_{i-1}$, the slice $M_{\mathbf p}(s_{i-1})$ has at most one neck above $a_i$. For $i=1$, this follows from the construction. For $i\ge2$, the profile of $M_{\mathbf p}(0)$ on $[a_i,a_{i-1}]$ agrees with that of $F_i$ and has at most one neck. The shrinking balls keep this interval positive for $0\le t\le2s_{i-1}$, and \eqref{eq:section-sign} keeps the derivatives at both endpoints negative. The Sturmian theorem therefore gives at most one neck there at $s_{i-1}$. The profile is strictly decreasing above $a_{i-1}$ at this time by the smooth convergence established above.

Since $u(a_i,t)$ stays positive and $u_y(a_i,t)<0$, no critical point crosses $y=a_i$. The neck bound consequently persists from $s_{i-1}$ until the first subsequent singular time or $2s_i$, whichever comes first. If this singular time occurs before $2s_i$, the component evolving from the connected slice $M_{\mathbf p}(s_{i-1})$ encloses the shrinking ball centered at the origin and cannot become extinct. The neck-loss argument of \cite{AAG}, restricted to $y>a_i$, shows that each distinct singular point at that time removes at least one neck above $a_i$.

\subsection{The Wa\.{z}ewski box argument}\label{sec:box}

For each $i$, consider the classical flow starting from $M_{\mathbf p}(s_{i-1})$, using the prescribed time $T_i$ and the constants fixed above. Parametrize this flow by
\[
 \ell=\log\frac{T_i-s_{i-1}}{T_i-t},\qquad \ell\ge0,
\]
and write $V_i(\mathbf p,\ell)$ for the coefficient vector in \eqref{eq:AV-low} along this flow. The box is defined by $|V_i(\mathbf p,\ell)|\le1$ for every $i$. For flows initially in the box, the continuation and estimates of Subsection~\ref{sec:AV-estimates} hold until the first exit, if any, and every finite first exit is strictly outward by \eqref{eq:AV-outward}.

\begin{proposition}\label{prop:box}
There exists $\mathbf p_*\in(B^3)^N$ such that, for every $i$, the classical flow from $M_{\mathbf p_*}(s_{i-1})$ exists on $[s_{i-1},T_i)$ and satisfies $|V_i(\mathbf p_*,\ell)|\le1$ for all $\ell\ge0$.
\end{proposition}
\begin{proof}
Let $\tau_i(\mathbf p)$ be the first $\ell\ge0$ such that $|V_i(\mathbf p,\ell)|\ge1$, and set $\tau_i=\infty$ if there is no such value. The continuation and strict outward crossing in Subsection~\ref{sec:AV-estimates}, including at the initial time, give continuity of $\tau_i$ with values in $[0,\infty]$ by classical continuous dependence, as in \cite[Lem.~3.2]{AV}. At $\tau_i=\infty$, continuity follows by applying continuous dependence on each compact $\ell$-interval.

Suppose every parameter makes at least one flow exit. Then $\tau(\mathbf p)=\min_i\tau_i(\mathbf p)$ is finite and continuous. Each flow exists up to $\ell=\tau(\mathbf p)$, and at least one of the vectors $V_i(\mathbf p,\tau(\mathbf p))$ has norm at least one. Define the exit map
\[
 \mathcal F(\mathbf p)=
 \frac{(V_1(\mathbf p,\tau(\mathbf p)),\ldots,V_N(\mathbf p,\tau(\mathbf p)))}
 {|(V_1(\mathbf p,\tau(\mathbf p)),\ldots,V_N(\mathbf p,\tau(\mathbf p)))|}
 \in\mathbb S^{3N-1}.
\]
This map is continuous, including at simultaneous exits. On the boundary of $(B^3)^N$, some $|\mathbf p_i|=1$, and \eqref{eq:AV-first-degree} and \eqref{eq:AV-second-degree} give
\[
 V_i(\mathbf p,0)\cdot\mathbf p_i>\tfrac32.
\]
Thus $\tau=0$ there, and the straight-line homotopy of the initial coefficient vector to $2\mathbf p$ never vanishes on the boundary. The restriction of $\mathcal F$ to the boundary is therefore homotopic to the radial map $\mathbf p\mapsto\mathbf p/|\mathbf p|$ and has degree one. This is impossible because $\mathcal F$ extends continuously to all of $(B^3)^N$ with values in the same sphere, so its boundary restriction must have degree zero.
\end{proof}

\begin{proof}[Proof of Theorem~\ref{thm:main}]
Choose $\mathbf p_*$ as in Proposition~\ref{prop:box} and write $M(t)=M_{\mathbf p_*}(t)$. For every $i$, uniqueness of the level set flow implies that $M(t)$ agrees on $[s_{i-1},T_i)$ with the classical flow starting from $M(s_{i-1})$. The argument for \cite[Thm.~3.3, (3.15)]{AV}, with the cutoff estimates in Subsection~\ref{sec:AV-estimates}, gives a limiting cusp at $(0,q_i)$ at time $T_i$, so this point is singular. The neck bound and the exclusion of singularities in $y\le a_i$ established in Subsection~\ref{sec:initial-surfaces}, together with the singularity description in \cite{AAG}, show that this is the only singular point at $T_i$.

The Angenent--Vel\'azquez estimates show that the part of the surface with $y\ge q_i$ shrinks to $(0,q_i)$, while the cusp profile and the shrinking ball centered at the origin ensure that a nonempty component survives below $q_i$. The cusp point $(0,q_i)$ lies in the closure of the surviving lower surface, so $M(T_i)$ is connected. By \cite{AAG}, the surviving component is smooth immediately after $T_i$; the neck count above $a_i$ shows that no neck remains there. Its upper cap is therefore monotone by \eqref{eq:section-sign}. The positive lower bound for $u(a_i,t)$ and the graphical interior estimates used for the lower cap in Subsection~\ref{sec:initial-surfaces} exclude further singularities above $a_i$ through $2s_i$. Together with the regularity in $y\le a_i$, this gives smooth connected evolution on $(T_i,2s_i]$. For $i<N$, this smooth interval overlaps the classical interval $[s_i,T_{i+1})$. Thus $T_1,\ldots,T_N$ are exactly the first $N$ singular times. Since there are initially at most $N$ necks and at least one is lost at each $T_i$, the component remaining after $T_N$ has no neck and stays smooth and connected until extinction at one point.

We express the asymptotics in the original coordinate $y$, using $\rho=(y-q_i)/\sqrt{T_i-t}$ and $\xi=(y-q_i)/(T_i-t)^{1/3}$. These coordinates are the negatives of those in Subsection~\ref{sec:AV-estimates}. For each $i$ there is a constant $K_i>0$ such that, as $t\uparrow T_i$,
\begin{align}
 \frac{u(q_i+\sqrt{T_i-t}\,\rho,t)}{\sqrt{T_i-t}}
 &=\sqrt2-\frac{K_i}{2\sqrt2}\sqrt{T_i-t}
       \operatorname{Hm}_3(\rho)+o(\sqrt{T_i-t}),
       \label{eq:AV-parabolic}\\
 \frac{u(q_i+(T_i-t)^{1/3}\xi,t)}{\sqrt{T_i-t}}
 &\longrightarrow\sqrt{2-K_i\xi^3}.
       \label{eq:AV-intermediate}
\end{align}
The first expansion holds uniformly on compact $\rho$-intervals, and the second convergence holds uniformly on compact subsets where $2-K_i\xi^3>0$. Set $\xi_i=(2/K_i)^{1/3}$, and write $H_i(t)$ for the mean curvature at the upper pole. The estimates in \cite[Sec.~7]{AV} give
\[
 b(t)=q_i+(\xi_i+o(1))(T_i-t)^{1/3},\qquad
 H_i(t)
   =\left(\frac{\xi_i}{3}+o(1)\right)(T_i-t)^{-2/3}.
\]
By \eqref{eq:AV-parabolic}, the parabolic limit is cylindrical. At the upper pole, which converges to $(0,q_i)$, we have $(T_i-t)|A|^2\ge\tfrac12(T_i-t)H_i(t)^2\to\infty$. Thus each singularity is locally Type~II and hence a degenerate neckpinch by \cite[Thm.~1.8]{SS}.
\end{proof}

\bibliographystyle{amsalpha}
\bibliography{main}
\end{document}